\documentclass[12pt,eqno]{article}

\usepackage{latexsym,amsmath,amssymb,amsthm,amsfonts}

\renewcommand{\epsilon}{\varepsilon}

\begin{document}

\author{Andriy V.\ Bondarenko and Maryna S. Viazovska}
\title{New asymptotic estimates for spherical designs}
\date{}
\maketitle \noindent Department of Mathematical Analysis, Kyiv
Taras Shevchenko University, Volodymyrska, 01033, Kyiv, Ukraine\\
tel: +38-044-259-05-91\\ fax:+38-044-259-03-92\\ e-mail:
bonda@univ.kiev.ua \\ Max Planck Institute for Mathematics,
Vivatsgasse 7, 53111 Bonn, Germany\\ tel:+49-228-402-265\\ e-mail:
viazovsk@mpim-bonn.mpg.de
\newpage
\begin{abstract}
Let $N(n,t)$ be the minimal number of points in a spherical
$t$-design on the unit sphere $S^n$ in $R^{{n+1}}$. For each $n\ge
3$, we prove a new asymptotic upper bound
$$
N(n,t)\le C(n)t^{a_n},
$$
where $C(n)$ is a constant depending only on $n$, $a_3\le 4$,
$a_4\le 7$, $a_5\le 9$, $a_6\le 11$, $a_7\le 12$, $a_8\le 16$,
$a_9\le 19$, $a_{10}\le 22$, and
$$
a_n< \frac n2\log_22n, \quad n> 10.
$$
\end{abstract}
\newpage
\section {Introduction}
Let $S^n$ be the unit sphere in $R^{{n+1}}$. The following concept
of a spherical design was introduced by Delsarte, Goethals and
Seidel ~\cite{DGS}:\\ The set of vectors
$\vec{x}_1,\ldots,\vec{x}_N\in S^n$ is called a {\sl spherical
$t$-design} if $$ \frac 1{mes\,
S^n}\int_{S^n}p(\vec{x})d\vec{x}=\frac 1N\sum_{i=1}^Np(\vec{x}_i)
$$ for all algebraic polynomials in $n+1$ variables and of total degree
$\le t$. For each $t\in {\mathbb N}$ denote by $N(n,t)$ the minimal
number of points in spherical $t$-design.  The following low bound
\begin{equation}
\label{hh}
N(n,t)\ge {{n+k}\choose{n}}+{{n+k-1}\choose{n}},\quad
t=2k,
\end{equation}
$$
N(n,t)\ge 2\,{{n+k}\choose{n}}, \quad t=2k+1,
$$
is also proved in \cite{DGS}. Spherical $t$-designs attaining these
bounds are called tight. Exactly eight tight spherical designs are
known for $n\ge 2$ and $t\ge 4$. All of them are obtained from
discrete algebraic structure called {\sl lattices}. Say, tight
7-design in $S^7$ is obtained from root $E_8$ lattice, and tight
11-design in $S^{23}$ is obtained from Leech lattice, see~\cite{CS}.
In general, lattices are a good source for spherical designs with
small $(n,t)$~\cite{PPV}. On the other hand construction of
spherical $t$-design with minimal cardinality for fixed $n$ and
$t\to\infty$ becomes a difficult analytic problem even for $n=2$.
There is strong relation between this problem and energy problem,
that is to find $N$ points on a sphere $S^2$ minimizing energy
functional $$ E(\vec{x}_1,\ldots,\vec{x}_N)=\sum_{1\le i<j\le
N}\frac 1{\|\vec{x}_i-\vec{x}_j\|},$$ see Saff, Kuijlaars~\cite{SK}.
Now we give a short history on asymptotic upper bounds on $N(n,t)$
for fixed $n$ and $t\to\infty$. First Seymour and Zaslavsky
\cite{SZ} have proved that spherical design exists for all $n$,$t\in
{\mathbb N}$. Then, Wagner \cite{W} and Bajnok \cite{B}
independently have proved that $N(n,t)\le C(n)t^{O(n^4)}$ and
$N(n,t)\le C(n)t^{O(n^3)}$ respectively. Korevaar and Meyers
\cite{KM} improved this inequality to $N(n,t)\le C(n)t^{(n^2+n)/2}$.
They also have conjectured that $N(n,t)\le C(n)t^n$. Remark,
\eqref{hh} implies $N(n,t)\ge C_1(n)t^n$. The main result of this
paper is \\ {\bf Theorem 1.} Let $a_n$ be the sequence defined by $$
a_1=1,\quad a_2=3, \quad a_{2n-1}=2a_{n-1}+n, \quad
a_{2n}=a_{n-1}+a_{n}+n+1,\quad n\ge 2. $$ Then, for all $n$,$t\in
{\mathbb N}$ we have
$$
N(n,t)\le C(n)t^{a_n},
$$
where $C(n)$ is a constant depending only on $n$.
\\ {\bf Corollary 1.} For each $n\ge 3$ and $t\in {\mathbb N}$ we have
$$
N(n,t)\le C(n)t^{a_n},
$$
where $C(n)$ is a constant depending only on $n$,
$$a_3\le 4, \quad a_4\le 7, \quad a_5\le 9, \quad a_6\le 11, \quad a_7\le 12, \quad a_8\le 16, \quad a_9\le 19, \quad a_{10}\le
22,$$ and
\begin{equation}
\label{hh1} a_n< \frac n2\log_22n, \quad n> 10.
\end{equation}

To prove Theorem 1 we need some auxiliary results.
\section{Auxiliary results}
We begin this section with following definition. Let $\omega(x)$ be
an integrable function on $[-1;1]$, then the set of points $T=
\{t_1,\ldots,t_K\}$, $t_k\in [-1,1]$, $k=\overline{1,K}$ is called
Chebyshev-type quadrature of degree $t$ with weight $\omega$ if $$
\int_{-1}^1p(x)\omega(x)dx=\frac 1K\int_{-1}^1\omega(x)dx
\sum_{k=1}^Kp(t_k) $$ for all algebraic polynomials $p$ of degree at
most $t$. An important class~of~weight functions are  Jacobi weight
functions $\omega_{m,n}=(1-x)^{(m-2)/2}(1+x)^{(n-2)/2}$. In
\cite{KK} Kuijlaars proved that for all $m$,$n\in {\mathbb N}$ there
exists Chebyshev-type quadrature of degree $t$ with weight
$\omega_{m,n}$ having at most $c(m,n)t^{\max(m,n)}$ points, where
$c(m,n)$ depends only on $m$ and $n$. For the general reference see
also \cite{DT}. Theorem 1 follows from this result and
\\ {\bf Lemma 1.} Let $X=\{\vec{x}_1,\ldots,\vec{x}_M\}$ and
$Y=\{\vec{y}_1,\ldots,\vec{y}_N\}$  be  spherical $t$-designs on
$S^{m-1}$ and $S^{n-1}$ respectively and $T= \{t_1,\ldots,t_K\}$ be
Chebyshev-type quadrature of degree $t$ with weight $\omega_{m,n}$.
Then there exists a spherical $t$-design on
$S^{m+n-1}$, having at most $KMN$ points.\\
Proof of Lemma 1 is based on the observation that any vector
$\vec{z}\in S^{m+n-1}$ can be written as
$$
\vec{z}=(\vec{x}\sin{\alpha}, \vec{y}\cos{\alpha}),
$$
where $\vec{x}\in S^{m-1}$, $\vec{y}\in S^{n-1}$ and $\alpha\in [0,
\pi/2]$. This decomposition is unique, except if $\vec{x}= 0$ or
$\vec{y} =0$, and it allows to transfer integration on $S^{n+m-1}$
to integration on the product space $$\Pi=S^{m-1}\times
S^{n-1}\times [0, \pi/2]$$ with Lebesgue measure on $S^{m-1}$ and a
measure on $[0, \pi/2]$ that can be mapped to the Jacobi measure
$\omega_{m,n}dx$ on [-1,1]. The main idea is the fact that these
transformations map even polynomials to polynomials of the same
degree on the product space. So, the corresponding spherical designs
on $S^{m-1}$ and $S^{n-1}$ and Chebyshev-type quadrature yield
spherical design on $S^{n+m-1}$.

\section{Proofs}
{\bf Proof of Lemma 1.} We will prove that a required spherical
$t$-design on $S^{n+m-1}$, say, is the set $L$ consisting of $KMN$
vectors of the form
$$ \vec{z}=\left(\sqrt{\frac{1-t}2}x_{1},\ldots,
\sqrt{\frac{1-t}2}x_{m},\sqrt{\frac{1+t}2}y_{1},\ldots,
\sqrt{\frac{1+t}2}y_{n}\right), $$ where $t \in T$,
$\vec{x}=(x_1,...,x_m)\in X$ and $\vec{y}=(y_1,...,y_n)\in Y$. To
this end take an arbitrary monomial $p(\vec{z})=z_1^{\alpha_1}\ldots
z_{m+n}^{\alpha_{m+n}}$ such that $\sum_{i=1}^{m+n}\alpha_i\le t$.
Note, that
$$
\sum_{\vec{z}\in
L}p(\vec{z})=\sum_{i=1}^M\sum_{j=1}^N\sum_{k=1}^K\prod_{l=1}^m
\left(x_{il}\sqrt{\frac{1-t_k}2}\right)^{\alpha_l}\prod_{s=1}^n
\left(y_{js}\sqrt{\frac{1+t_k}2}\right)^{\alpha_{m+s}}
$$
$$
=\sum_{i=1}^M\prod_{l=1}^mx_{il}^{\alpha_l}
\sum_{j=1}^N\prod_{s=1}^ny_{js}^{\alpha_{m+s}}
\sum_{k=1}^K\left(\sqrt{\frac{1-t_k}2}\right)
^{\alpha_1+\ldots+\alpha_m}\left(\sqrt{\frac{1+t_k}2}\right)
^{\alpha_{m+1}+\ldots+\alpha_{m+n}}.
$$
Set
\begin{equation}
\label{2} I_1=\frac
1{mes\,S^{m-1}}\int_{S^{m-1}}x_1^{\alpha_1}\ldots
x_{m}^{\alpha_{m}}dx_1\ldots dx_m,
\end{equation}
\begin{equation}
\label{3} I_2=\frac
1{mes\,S^{n-1}}\int_{S^{n-1}}y_{1}^{\alpha_{m+1}}\ldots
y_{n}^{\alpha_{m+n}}dy_{1}\ldots dy_{n},
\end{equation}
and $$ I_3=\frac 1K\sum_{k=1}^K\left(\sqrt{\frac{1-t_k}2}\right)
^{\alpha_1+\ldots+\alpha_m}\left(\sqrt{\frac{1+t_k}2}\right)
^{\alpha_{m+1}+\ldots+\alpha_{m+n}}. $$ Since
$\vec{x}_1,\ldots,\vec{x}_N$ and $\vec{y}_1,\ldots,\vec{y}_M$ are
spherical $t$-designs, we have
\begin{equation}
\label{1} \frac 1{KMN}\sum_{\vec{z}\in L}p(\vec{z})=I_1I_2I_3,
\end{equation}
If, for some $i=\overline{1,m+n}$, $\alpha_i$ is odd, then either
$I_1=0$ or $I_2=0$, hence $$ \frac 1{KMN}\sum_{\vec{z}\in
L}p(\vec{z})=\frac
1{mes\,S^{m+n-1}}\int_{S^{m+n-1}}p(\vec{z})d\vec{z}=0. $$ So, we may
assume that all $\alpha_i$ are even. Put $\beta_i:=\alpha_i/2$,
$i=\overline{1,m+n}$. Since $t_1,\ldots,t_K\in [-1,1]$ is a
Chebyshev-type quadrature of degree $t$ with weight $\omega_{m,n}$,
then $$ I_3=\frac 1K\sum_{k=1}^K\left(\frac{1-t_k}2\right)
^{\beta_1+\ldots+\beta_m}\left(\frac{1+t_k}2\right)
^{\beta_{m+1}+\ldots+\beta_{m+n}} $$ $$ =\frac
1{\int_{-1}^1\omega_{m,n}(t)dt}\int_{-1}^1\left(\frac{1-t}2\right)
^{\beta_1+\ldots+\beta_m}\left(\frac{1+t}2\right)
^{\beta_{m+1}+\ldots+\beta_{m+n}}\omega_{m,n}(t)dt. $$ Using change
of variables $t=\cos{2\alpha}$, $\alpha\in [0,\pi/2]$, we obtain
\begin{equation}
\label{4} I_3=\frac
{2^{\frac{m+n}{2}}}{\int_{-1}^1\omega_{m,n}(t)dt}\int_0^{\pi/2}(\sin{\alpha})
^{\alpha_1+\ldots+\alpha_m+m-1}(\cos{\alpha})^{\alpha_{m+1}+\ldots+\alpha_{m+n}+n-1}d\alpha.
\end{equation}
Now we are ready to prove the equality
\begin{equation}
\label{12} \frac
1{mes\,S^{m+n-1}}\int_{S^{m+n-1}}p(\vec{z})d\vec{z}=\frac
1{KMN}\sum_{\vec{z}\in L}p(\vec{z}).
\end{equation}
To this end we introduce change of variables
$$
\vec{z}=(\vec{x}\sin{\alpha}, \vec{y}\cos{\alpha}),
$$
where
$$(\vec{x},\vec{y},\alpha)\in S^{m-1}\times S^{n-1}\times [\,0,
\pi/2\,]=:\Pi.$$ Jacobian $J$ of this transformation is
$$
J(\vec{x},\vec{y},\alpha)=(\sin{\alpha})^{m-1}(\cos{\alpha})^{n-1},
$$
therefore, by \eqref{2}-\eqref{4},
$$ \int_{S^{m+n-1}}p(\vec{z})d\vec{z}=
\int_{\Pi}p((\vec{x}\sin{\alpha},\vec{y}\cos{\alpha}))(\sin{\alpha})^{m-1}(\cos{\alpha})^{n-1}d\vec{x}d\vec{y}d\alpha
$$
$$
=\int_{\Pi}(x_1\sin{\alpha})^{{\alpha}_1}\ldots
(x_m\sin{\alpha})^{{\alpha}_m}(y_1\cos{\alpha})^{{\alpha}_{m+1}}
\ldots (y_n\cos{\alpha})^{{\alpha}_{m+n}}
(\sin{\alpha})^{m-1}(\cos{\alpha})^{n-1}d\vec{x}d\vec{y}d\alpha
$$
$$
=\int_{S^{m-1}}x_1^{\alpha_1}\ldots x_{m}^{\alpha_{m}}dx_1\ldots
dx_m \int_{S^{n-1}}y_{1}^{\alpha_{m+1}}\ldots
y_{n}^{\alpha_{m+n}}dy_{1}\ldots dy_{n}
$$
$$
 \int_0^{\pi/2}(\sin{\alpha})
^{\alpha_1+\ldots+\alpha_m+m-1}(\cos{\alpha})^{\alpha_{m+1}+\ldots+\alpha_{m+n}+n-1}d\alpha
= A_{m,n}\,I_1I_2I_3,
$$
where
$$ A_{m,n}=
2^{-\frac{m+n}2}\,{mes\,S^{m-1}}{mes\,S^{n-1}}\int_{-1}^1\omega_{m,n}(t)dt.
$$ Since \eqref{12} holds for $p(\vec{z})=1$, then $A_{m,n}=mes\,S^{m+n-1}$ for
all $m,n\in {\mathbb N}$, hence $L$ is a spherical $t$-design on
$S^{m+n-1}$. Lemma 1 is proved.\\ {\bf Proof of Theorem 1.} We prove
Theorem 1 by induction on n.
\\ If $n=1$, then $N(1,t)=t+1\le C(1)t$. For $n=2$ the estimate
$N(2,t)\le C(2)t^3$ is proved in \cite{KM}. Lemma 1 and existence of
Chebyshev-type quadrature of degree $t$ with weight $\omega_{m,n}$
having at most $c(m,n)t^{\max(m,n)}$ points imply
$$N(n+m-1,t)\leq N(n-1,t)N(m-1,t)c(n,m)t^{\max(m,n)}.$$ So, taking
either $m=n$ or $m=n+1$, we get $$N(2n-1,t)\leq
C(n-1)C(n-1)c(n,n)t^{2a_{n-1}+n}=C^2(n-1)c(n,n)t^{a_{2n-1}}=:C(2n-1)t^{a_{2n-1}}$$
and
$$N(2n,t)\leq
C(n-1)C(n)c(n-1,n)t^{a_{n-1}+a_n+n+1}=C(n-1)C(n)c(n-1,n)t^{a_{2n}}=:C(2n)t^{a_{2n}}.$$
Theorem 1 is proved.\\ {\bf Proof of Corollary 1.} For $n\le 21$ one
check \eqref{hh1} directly. For $n> 21$ one checks \eqref{hh1} by
induction.\\
{\bf Remark.} By the definition of the sequence
$\{a_n\}_{n=1}^{\infty}$,
$$ a_{2^{n}-1}=n2^{n-1}, \quad n\in{\mathbb N},$$ hence
$$ \limsup_{n\to\infty}\frac{a_n}{n\log_2n}=\frac 12. $$ So, we
cannot improve the constant $1/2$ in \eqref{hh1}. Corollary 1 is proved.\\
{\bf Acknowledgement.} The authors would like to thank Professor
Edward Saff for his inspiring lecture on the energy problems, Kyiv,
June 2006. We are most grateful to Professor Arno Kuijlaars for
fruitful remarks that allowed to improve the text essentially, and,
moreover, to strengthen the main result.

\end{document}